# PREHOMOGENEOUS VECTOR SPACES AND ERGODIC THEORY II

DAVE WITTE[1], AKIHIKO YUKIE[2], AND ROGER ZIERAU[3]

Oklahoma State University

## Contents

## Introduction

Throughout this paper, $k$ is a field of characteristic zero. In this paper we consider the following three prehomogeneous vector spaces
(1) $G = \text{GL}(6)$, $V = \wedge^3 k^6$,
(2) $G = \text{GL}(1) \times \text{GL}(7)$, $V = \wedge^3 k^7$,
(3) $G = \text{GL}(2n)$, $V = \wedge^2 k^{2n}$.
For the definition of prehomogeneous vector spaces, see [7], [11], or [10]. The classification of irreducible regular prehomogeneous vector spaces is known [7].

For any algebraic group $G$, we denote the connected component of 1 in Zariski topology by $G^0$. If $G$ is defined over a subfield of $\mathbb{R}$, we denote the connected component of 1 of $G_\mathbb{R}$ in classical topology by $G_{\mathbb{R}+}^0$. Let $(G, V)$ be one of the prehomogeneous vector spaces (1)–(3). Let $H = \text{SL}(n) \subset G$ ($n = 6, 7$) or $H = \text{SL}(2n)$. For $x \in V_k^{\text{ss}}$, let $G_x$ be the stabilizer and $H_x = G_x \cap H$.

Consider the case $k = \mathbb{Q}$. Note that $H_{\mathbb{R}+}^0 = H_\mathbb{R}$. Let $\Gamma \subset H_{\mathbb{R}+}^0 = H_\mathbb{R}$ be an arithmetic lattice. The second author pointed out in [10], as a consequence of the Moore ergodicity theorem, that for almost all $x \in V_\mathbb{R}^{\text{ss}}$ such that $H_{x\mathbb{R}}^0$ has positive real rank, the set $H_{x\mathbb{R}+}^0 \Gamma$ is dense in $H_\mathbb{R}$. The purpose of this paper is to find an explicit irrationality condition on $x$ that implies the density of this set, and to

---

*Corresponding author:* A. Yukie
[1] Partially supported by NSF grant DMS-9214077
[2] Partially supported by NSF grant DMS-9401391
[3] Partially supported by NSF grant DMS-9303224

Typeset by $\mathcal{AMS}$-TEX

interpret the result number theoretically. This provides an answer to part (2) of Question (0.5) in [10], for the three prehomogeneous vector spaces that we consider.

The following theorem is the main result of this paper. In statements (1)–(3) of the following theorem, we consider the prehomogeneous vector spaces (1)–(3) respectively. Let $W$ be the standard representation of $\mathrm{GL}(6), \mathrm{GL}(7)$, or $\mathrm{GL}(2n)$. We identify $\wedge^3 W$ and $\wedge^2 W$ with $\wedge^3(W^*)^*$ and $\wedge^2(W^*)^*$ respectively. We define the notion of "sufficiently irrational" points in Definition (5.1)(1)–(3).

**Theorem (6.15)** *Suppose $x \in V_{\mathbb{R}}^{\mathrm{ss}}$ is sufficiently irrational and the real rank of $H_{x\mathbb{R}}^0$ is positive.*
*(1) For any $y = (y_{ijk}) \in \wedge^3 \mathbb{R}^5$ and $\epsilon > 0$, there exists a $\mathbb{Z}$–basis $\{u_1, \cdots, u_6\}$ of $W_{\mathbb{Z}}^* \cong \mathbb{Z}^6$ such that*
$$|y_{ijk} - x(u_i, u_j, u_k)| < \epsilon$$
*for all $i < j < k \leq 5$.*
*(2) For any $y = (y_{ijk}) \in \wedge^3 \mathbb{R}^6$ and $\epsilon > 0$, there exists a $\mathbb{Z}$–basis $\{u_1, \cdots, u_7\}$ of $W_{\mathbb{Z}}^* \cong \mathbb{Z}^7$ such that*
$$|y_{ijk} - x(u_i, u_j, u_k)| < \epsilon$$
*for all $i < j < k \leq 6$.*
*(3) For any $y = (y_{ij}) \in \wedge^2 \mathbb{R}^{2n-1}$ and $\epsilon > 0$, there exists a $\mathbb{Z}$–basis $\{u_1, \cdots, u_{2n}\}$ of $W_{\mathbb{Z}}^* \cong \mathbb{Z}^{2n}$ such that*
$$|y_{ij} - x(u_i, u_j)| < \epsilon$$
*for all $i < j \leq 2n - 1$.*

## §1 The orbit space $G_k \backslash V_k^{\mathrm{ss}}$ (1)

We describe the orbit space $G_k \backslash V_k^{\mathrm{ss}}$ in §§1,2. We do not need the GL(1)–factor for cases (1) and (3), because having an extra GL(1)–factor does not change the orbit space $G_k \backslash V_k^{\mathrm{ss}}$ for these cases. Let $W$ be as in the introduction. Throughout this paper, for a basis $\{e_1, e_2, \cdots\}$ of $W$, we use the notation $e_{i_1 \cdots i_k} = e_{i_1} \wedge \cdots \wedge e_{i_k}$.

For case (3), $G_k \backslash V_k^{\mathrm{ss}}$ consists of a single point due to the well known fact that over any field two symplectic forms are equivalent. Let $e_1, \cdots, e_{2n}$ be the standard coordinate vectors of $k^{2n}$. Let $w = e_{1\ n+1} + \cdots + e_{n\ 2n}$. Then $G_k \backslash V_k^{\mathrm{ss}} = G_k w$ and $G_w = G_w^0 = H_w = H_w^0 = \mathrm{Sp}(2n)$.

We consider case (1) for the rest of this section. Let $G = \mathrm{GL}(6)$, $W = k^6$, and $V = \wedge^3 W$. It is known (see [7, p. 80]) that this is a prehomogeneous vector space. Let $\{e_1, \cdots, e_6\}$ be a basis of $W$. It is known (see [7]) that the orbit of

(1.1) $$w = e_{123} + e_{456}$$

is Zariski open in $V$.

Let

(1.2) $$d(A, B) = \begin{pmatrix} A & 0 \\ 0 & B \end{pmatrix}$$

for $A, B \in \mathrm{M}(3,3)$. Then

(1.3) $$G_w^0 = \{d(A, B) \mid A, B \in \mathrm{GL}(3),\ \det A = \det B = 1\} \cong \mathrm{SL}(3) \times \mathrm{SL}(3)$$

and $G_w$ is generated by $G_w^0$ and the element

$$\tag{1.4} \tau = \begin{pmatrix} 0 & I_3 \\ I_3 & 0 \end{pmatrix}.$$

This implies that the following sequence

$$\tag{1.5} 1 \to G_w^0 \to G_w \to \mathbb{Z}/2\mathbb{Z} \to 1$$

which sends $\tau$ to the non-trivial element of $\mathbb{Z}/2\mathbb{Z}$ is exact. Moreover, this exact sequence is split. For these facts, the reader should see [7]. Note that the argument in [7] works over any ground field of characteristic zero.

For any algebraic group $G$ over $k$, let $\mathrm{H}^1(k, G)$ be the first Galois cohomology set. We choose the definition so that trivial classes are those of the form $\{g^{-1}g^\sigma\}_{\sigma \in \mathrm{Gal}(\bar{k}/k)}$ ($g \in G_{\bar{k}}$) and the cocycle condition is $h_{\sigma\tau} = h_\tau h_\sigma^\tau$ for a continuous map $\{h_\sigma\}_{\sigma \in \mathrm{Gal}(\bar{k}/k)}$ from $\mathrm{Gal}(\bar{k}/k)$ to $G_{\bar{k}}$.

For $x = gw \in V_k^{\mathrm{ss}}$ where $g \in G_{\bar{k}}$, we define $c_x = \{g^{-1}g^\sigma\}_{\sigma \in \mathrm{Gal}(\bar{k}/k)} \in \mathrm{H}^1(k, G_w)$. This definition does not depend on the choice of $g$. Since $\mathrm{H}^1(k, G) = \{1\}$, by [4, p. 269],

$$\tag{1.6} G_k \backslash V_k^{\mathrm{ss}} \ni x \to c_x \in \mathrm{H}^1(k, G_w)$$

is bijective.

Let $\mathfrak{Ex}_2$ be the set of isomorphism classes of extensions $k'/k$ of degree either one or two. By the split exact sequence (1.5), we get a surjective map

$$\tag{1.7} \alpha_V : G_k \backslash V_k^{\mathrm{ss}} \to \mathrm{H}^1(k, \mathbb{Z}/2\mathbb{Z}) \cong \mathfrak{Ex}_2.$$

For $x \in G_k \backslash V_k^{\mathrm{ss}}$, we denote the field corresponding to $\alpha_V(x)$ by $k(x)$. Let $k(\alpha)$ be the field generated by an element of the form $\alpha = \sqrt{\beta}$. We define $w_\alpha = g_\alpha w$ where

$$\tag{1.8} g_\alpha = \begin{pmatrix} 1 & 0 & 0 & 1 & 0 & 0 \\ \alpha & 0 & 0 & -\alpha & 0 & 0 \\ 0 & 1 & 0 & 0 & 1 & 0 \\ 0 & \alpha & 0 & 0 & -\alpha & 0 \\ 0 & 0 & 1 & 0 & 0 & 1 \\ 0 & 0 & \alpha & 0 & 0 & -\alpha \end{pmatrix}.$$

Suppose $\alpha \notin k^\times$ and $\sigma$ is the non-trivial element of $\mathrm{Gal}(k(\alpha)/k)$. It is easy to see that $g_\alpha^\sigma = g_\alpha \tau$. This implies that $(g_\alpha w)^\sigma = g_\alpha \tau w = g_\alpha w$. So $w_\alpha \in V_k^{\mathrm{ss}}$ and $k(w_\alpha) = k(\alpha)$. Easy computations show that

$$\tag{1.9} w_\alpha = e_{123} + \alpha^2(e_{156} - e_{246} + e_{345}).$$

We determine $G_{w_\alpha}^0$ rationally.

**Proposition (1.10)** If $\alpha \notin k^\times$, as an algebraic group over $k$,

$$G_{w_\alpha}^0 = \{g_\alpha d(A, A^\sigma) g_\alpha^{-1} \,|\, A \in \mathrm{SL}(3)_{k(\alpha)}\}$$
$$\cong \mathrm{SL}(3)_{k(\alpha)}.$$

*Proof.* In order to prove an isomorphism between two algebraic groups $G_1, G_2$ over $k$, it is enough to prove natural isomorphisms between the sets $G_{1R}, G_{2R}$ of $R$–rational points of $G_1, G_2$ for all $k$–algebras $R$. For this, the reader should see Theorem [6, p. 17].

Let $R$ be any $k$–algebra. For any finite Galois extension $k'/k$, $\nu \in \mathrm{Gal}(k'/k)$ acts on $k' \otimes R$ by $(x \otimes r)^\nu = x^\nu \otimes r$. We define $R(\alpha) = k(\alpha) \otimes R$. Then
$$G^0_{w_\alpha R} = \{g \in G^0_{w_\alpha R(\alpha)} \mid g^\sigma = g\}.$$

Over $R(\alpha)$, we can express elements of $G^0_{w_\alpha R}$ as

(1.11) $$g_\alpha d(A, B) g_\alpha^{-1},$$

where $A, B \in \mathrm{SL}(3)_{R(\alpha)}$. The element (1.11) is in $G^0_{w_\alpha R}$ if and only if
$$\begin{aligned} g_\alpha d(A, B) g_\alpha^{-1} &= g_\alpha^\sigma d(A^\sigma, B^\sigma)(g_\alpha^\sigma)^{-1} \\ &= g_\alpha \tau d(A^\sigma, B^\sigma) \tau g_\alpha^{-1} \\ &= g_\alpha d(B^\sigma, A^\sigma) g_\alpha^{-1}. \end{aligned}$$

This condition is satisfied if and only if $B = A^\sigma$. This proves the proposition. $\square$

Note that
$$g_\alpha \tau g_\alpha^{-1} = \begin{pmatrix} I_3 & \\ & -I_3 \end{pmatrix}.$$

So the sequence
$$1 \to G^0_{w_\alpha} \to G_{w_\alpha} \to \mathbb{Z}/2\mathbb{Z} \to 1$$
is also a split exact sequence.

By Lemma (1.4) [5], $\mathrm{H}^1(k, G^0_{w_\alpha})$ is trivial for all $\alpha$. So by Lemma (1.12) [5], we get the following proposition.

**Proposition (1.12)** *The map $\alpha_V$ induces a bijection $G_k \backslash V_k^{\mathrm{ss}} \cong \mathfrak{Er}_2$.*

This gives us an interpretation of the expected density theorem from the zeta function theory of this case and the zeta function is a counting function of
$$\operatorname*{Res}_{s=1} \zeta_{k(\alpha)}(s) \zeta_{k(\alpha)}(2) \zeta_{k(\alpha)}(3)$$
for quadratic extensions $k(\alpha)$ ($\zeta_{k(\alpha)}(s)$ is the Dedekind zeta function). However, we will not consider the zeta function in this paper. For the zeta function theory of prehomogeneous vector spaces, see [9], [11].

If $k = \mathbb{R}$, $G_\mathbb{R} \backslash V_\mathbb{R}^{\mathrm{ss}}$ is represented by $w_0 = w$ and

(1.13) $$w_1 = e_{123} - e_{156} + e_{246} - e_{345}.$$

As in [7], we define a map $D_3 : \wedge^3 W \to \wedge^2 W \otimes W$ by

(1.14) $$D_3(v_1 \wedge v_2 \wedge v_3) = v_2 \wedge v_3 \otimes v_1 - v_1 \wedge v_3 \otimes v_2 + v_1 \wedge v_2 \otimes v_3$$

for $v_1, v_2, v_3 \in W$.

For $x \in V_k$, we define

(1.15) $\quad S_x = x \wedge D_3(x) \in \wedge^5 W \otimes W \cong W^* \otimes W \cong \text{Hom}(W, W),$

where $x \wedge (y \otimes z) = (x \wedge y) \otimes z$.

It is proved in [7, p. 80] that there is a relative invariant polynomial $\Delta(x)$ of degree four such that $S_x^2 = \Delta(x) I_W$ where $I_W$ is the identity map of $W$. This implies that

(1.16) $\quad\quad\quad\quad\quad\quad \Delta(gx) = (\det g)^2 \Delta(x).$

Since $S_w$ has eigenvalues $\pm 1$ (see [7, p. 80]), $\Delta(w) = 1$. Also since $\operatorname{ch} k = 0$ and $\Delta(w) \in k$, we may assume that $\Delta \in k[V]$. Note that $\det g_\alpha = -8\alpha$. So $\Delta(w_\alpha) = 64\alpha^2 \Delta(w) = 64\alpha^2$. Therefore, this gives us the following characterization of the field $k(x)$.

**Proposition (1.17)** *For $x \in V_k^{\text{ss}}$, the field $k(x)$ is generated by eigenvalues of $S_x$.*

Let $E_{x1}, E_{x2}$ be the eigenspaces of $S_x$ for the eigenvalues $\pm \Delta(x)^{\frac{1}{2}}$. It is known [7] that $\dim E_{x1} = \dim E_{x2} = 3$. Let $\text{Gr}(3,6)$ be the Grassmann of 3–dimensional subspaces of $W$. Let

(1.18) $\quad\quad\quad\quad\quad X = (\mathbb{Z}/2\mathbb{Z}) \setminus (\text{Gr}(3,6) \times \text{Gr}(3,6)),$

where $\mathbb{Z}/2\mathbb{Z}$ acts by permuting two factors.

**Definition (1.19)** $\text{Gr}(x) = ([E_{x1}], [E_{x2}]) \in X$.

Since we are assuming $\operatorname{ch} k = 0$, $k$–rational points of $X$ are points which are set theoretically fixed by any $\sigma \in \text{Gal}(\bar{k}/k)$. If $x \in V_k^{\text{ss}}$, $E_{x1}, E_{x2}$ are conjugate. So $\text{Gr}(x)$ is a $k$–rational point of $X$. The following proposition is obvious.

**Proposition (1.20)** *If $x \in V_{\bar{k}}^{\text{ss}}$ and $\sigma \in \text{Gal}(\bar{k}/k)$, $\text{Gr}(x^\sigma) = \text{Gr}(x)^\sigma$.*

For general $x \in V^{\text{ss}}$, we cannot distinguish $E_{x1}$ and $E_{x2}$. But for later purposes, we choose $E_{w1}$ and $E_{w2}$ so that $E_{w1}$ is spanned by $e_1, e_2, e_3$ and $E_{w2}$ is spanned by $e_4, e_5, e_6$. It is proved in [7, p. 80] that $E_{w1}$ (resp. $E_{w2}$) is the eigenspace of $S_w$ for the eigenvalue $1$ (resp. $-1$).

## §2 The orbit space $G_k \setminus V_k^{\text{ss}}$ (2)

We describe the orbit space $G_k \setminus V_k^{\text{ss}}$ for case (2) in this section. Since this case has something to do with octonion algebras, we briefly recall the Cayley–Dickson process. Octonian algebras are often referred to as Cayley algebras also. For a reference, see [2, pp. 101–110] for example. Note that even though [2] assumes $k = \mathbb{R}$, the argument for the Cayley–Dickson process works over any field of characteristic zero.

**Definition (2.1)** *A normed $k$–algebra is a not necessarily associative finite dimensional $k$–algebra $A$ with multiplicative unit $1$, equipped with a non-degenerate* *symmetric bilinear form $\langle x, y \rangle$ for $x, y \in A$ such that the associated square norm $\|x\| = \langle x, x \rangle$ satisfies the multiplicative property*

$$\|xy\| = \|x\|\|y\|.$$

If $A$ is a normed $k$–algebra, we denote the span of 1 by $\mathrm{Re}(A)$ and its orthogonal complement $\{x \in A \mid \langle 1, x \rangle = 0\}$ by $\mathrm{Im}(A)$. Any $x \in A$ has a unique decomposition $x = x_1 + x_2$ with $x_1 \in \mathrm{Re}(A), x_2 \in \mathrm{Im}(A)$. We denote $\mathrm{Re}(x) = x_1, \mathrm{Im}(x) = x_2$. We define the conjugation by $\bar{x} = x_1 - x_2$. So $\mathrm{Re}(x) = \frac{1}{2}(x + \bar{x})$.

Given a normed $k$–algebra $A$, we make two new normed $k$–algebras $A(\pm)$ as follows. As a vector space,
$$A(\pm) = A \oplus A.$$
We define the multiplication and the norm by

(2.2)
$$(a, b)(c, d) = (ac \mp \bar{d}b, da + b\bar{c}),$$
$$\|(a, b)\| = \|a\| \pm \|b\|.$$

Then we define
$$\langle (a, b), (c, d) \rangle = \frac{1}{2}(\|(a, b) + (c, d)\| - \|(a, b)\| - \|(c, d)\|).$$

The algebra $A(\pm)$ becomes a normed $k$–algebra by the above product and the bilinear form. We use the notation $a + b\epsilon$ for $(a, b)$. Note that if $k$ contains $\sqrt{-1}$, $\epsilon \to \sqrt{-1}\epsilon$ induces an isomorphism $A(+) \to A(-)$.

For a normed $k$–algebra $A$, we define $[x, y, z] = (xy)z - x(yz)$ for $x, y, z \in A$. This is called the associator. If the associator is alternating, $A$ is called an alternative algebra. It is known that if $A$ is commutative, $A(\pm)$ is associative, and if $A$ is associative, $A(\pm)$ is alternative. The above process is called the Cayley–Dickson process. It is easy to see that

(2.3)
$$\mathrm{Im}(a + b\epsilon) = \mathrm{Im}(a) - b\epsilon, \; \overline{a + b\epsilon} = \bar{a} - b\epsilon.$$

The following lemma is proved in [2].

**Lemma (2.4)** (1) $\overline{xy} = \bar{y}\bar{x}$,
(2) $\langle x, y \rangle = \mathrm{Re}(x\bar{y})$.
(3) $\|x\| = x\bar{x}$.

If $A, B$ are normed $k$–algebras, a homomorphism $\phi : A \to B$ is a $k$–linear map such that $\phi(1) = 1$, $\phi(xy) = \phi(x)\phi(y)$, and $\|\phi(x)\| = \|x\|$. The third condition implies $\langle \phi(x), \phi(y) \rangle = \langle x, y \rangle$. So $\phi(\mathrm{Im}(A)) \subset \mathrm{Im}(B)$. Suppose $x, y \in \mathrm{Im}(A)$. Then $\langle x, y \rangle = \mathrm{Re}(x\bar{y}) = -\mathrm{Re}(xy)$. So

$$-\mathrm{Re}(\phi(x)\phi(y)) = \langle \phi(x), \phi(y) \rangle = \langle x, y \rangle = -\mathrm{Re}(xy),$$
$$\phi(\mathrm{Im}(xy)) = \phi(xy - \mathrm{Re}(xy)) = \phi(x)\phi(y) - \mathrm{Re}(xy)$$
$$= \phi(x)\phi(y) - \mathrm{Re}(\phi(x)\phi(y)) = \mathrm{Im}(\phi(x)\phi(y)).$$

Conversely, $\phi$ is a homomorphism if the above conditions are satisfied. So we have proved the following proposition.

**Proposition (2.5)** *A $k$–linear map $\phi : A \to B$ is a homomorphism if and only if*

$$\phi(1) = 1, \ \phi(\mathrm{Im}(xy)) = \mathrm{Im}(\phi(x)\phi(y)), \ \langle\phi(x), \phi(y)\rangle = \langle x, y\rangle$$

for all $x, y \in \mathrm{Im}(A)$.

It is easy to see that

(2.6)
$$k(+) \cong \left\{ \begin{pmatrix} a & -b \\ b & a \end{pmatrix} \middle| a, b \in k \right\},$$
$$k(+)(-) \cong \mathrm{M}(2,2)_k.$$

For $k(+)$, $\epsilon = \begin{pmatrix} 0 & -1 \\ 1 & 0 \end{pmatrix}$, and the conjugation is

$$\begin{pmatrix} a & -b \\ b & a \end{pmatrix} \to \begin{pmatrix} a & b \\ -b & a \end{pmatrix}.$$

For $k(+)(-)$, $\epsilon = \begin{pmatrix} 1 & 0 \\ 0 & -1 \end{pmatrix}$, and the conjugation is

$$\begin{pmatrix} a & b \\ c & d \end{pmatrix} \to \begin{pmatrix} d & -b \\ -c & a \end{pmatrix}.$$

Therefore, $\mathrm{Re}(x) = \frac{1}{2}\mathrm{tr}(x)$ in both cases and the norm is the determinant.

We define $\mathbb{H} = k(+)(+)$, $\mathbb{O} = \mathbb{H}(+)$, and $\widetilde{\mathbb{O}} = \mathrm{M}(2,2)(+)$. $\mathbb{O}$ is called the non-split octonion algebra (if $k$ does not contain $\sqrt{-1}$), and $\widetilde{\mathbb{O}}$ is called the split octonion algebra.

We consider properties of the split octonion algebra $\widetilde{\mathbb{O}}$. Let

(2.7)
$$C(x, y, z) = \langle x, yz \rangle$$

for $x, y, z \in \widetilde{\mathbb{O}}$. Suppose $x = x_1 + x_2\epsilon$, $y = y_1 + y_2\epsilon$, and $z = z_1 + z_2\epsilon$.

**Lemma (2.8)**

$$C(x, y, z) = \frac{1}{2}\mathrm{tr}(x_1(\bar{z}_1\bar{y}_1 - \bar{y}_2 z_2) + (\bar{y}_1 \bar{z}_2 + z_1 \bar{y}_2)x_2).$$

*Proof.* Since $yz = y_1 z_1 - \bar{z}_2 y_2 + (z_2 y_1 + y_2 \bar{z}_1)\epsilon$,

$$\overline{yz} = \bar{z}_1 \bar{y}_1 - \bar{y}_2 z_2 - (z_2 y_1 + y_2 \bar{z}_1)\epsilon,$$
$$x(\overline{yz}) = x_1(\bar{z}_1 \bar{y}_1 - \bar{y}_2 z_2) + (\bar{y}_1 \bar{z}_2 + z_1 \bar{y}_2)x_2$$
$$\quad + (-(z_2 y_1 + y_2 \bar{z}_1)x_1 + x_2(y_1 z_1 - \bar{z}_2 y_2))\epsilon.$$

Therefore,
$$\begin{aligned}
C(x,y,z) &= \mathrm{Re}(x(\overline{yz})) \\
&= \mathrm{Re}(x_1(\bar{z}_1\bar{y}_1 - \bar{y}_2 z_2) + (\bar{y}_1\bar{z}_2 + z_1\bar{y}_2)x_2) \\
&= \frac{1}{2}\mathrm{tr}(x_1(\bar{z}_1\bar{y}_1 - \bar{y}_2 z_2) + (\bar{y}_1\bar{z}_2 + z_1\bar{y}_2)x_2).
\end{aligned}$$

□

Note that $\mathrm{tr}(x) = \mathrm{tr}(\bar{x})$ for $x \in \mathrm{M}(2,2)$. Also if $x,y,z \in \mathrm{Im}(\widetilde{\mathbb{O}})$, $x_1 + \bar{x}_1 = y_1 + \bar{y}_1 = z_1 + \bar{z}_1 = 0$. Therefore, the above lemma and a straightforward argument shows that $C(x,y,z)$ is an element of $\wedge^3 \mathrm{Im}(\widetilde{\mathbb{O}})^*$ (see (6.60) [2, p. 113] also).

Let $E_{ij}$ be the $2 \times 2$ matrix whose $(i,j)$–entry is 1 and other entries are zero. Let

(2.9) $$f_1 = \begin{pmatrix} 1 & 0 \\ 0 & -1 \end{pmatrix}, \; f_2 = E_{12}, \; f_3 = E_{11}\epsilon, \; f_4 = -E_{21}\epsilon,$$
$$f_5 = -E_{21}, \; f_6 = E_{22}\epsilon, \; f_7 = E_{12}\epsilon.$$

Then $\{f_1, \cdots, f_7\}$ is a basis for $\mathrm{Im}(\widetilde{\mathbb{O}})$. Let $e_1, \cdots, e_7$ be the dual basis. Straightforward computations (35 computations for (2)) using the above lemma shows the following proposition and the proof is left to the reader.

**Proposition (2.10)** (1) *Suppose*

$$x = \sum_{i=1}^{7} x_i f_i = \begin{pmatrix} x_1 & x_2 \\ -x_5 & -x_1 \end{pmatrix} + \begin{pmatrix} x_3 & x_7 \\ -x_4 & x_6 \end{pmatrix}\epsilon.$$

*Then $\|x\| = -x_1^2 + x_2 x_5 + x_3 x_6 + x_4 x_7$.*
*(2) $C = \frac{1}{2}(e_{234} + e_{567} + e_{125} + e_{136} + e_{147})$.*

Let $W = \mathrm{Im}(\widetilde{\mathbb{O}})^*$. We choose the above basis $\{e_1, \cdots, e_7\}$. Let $G = \mathrm{GL}(1) \times \mathrm{GL}(W)$, $V = \wedge^3 W$. The GL(1)–factor acts by the usual multiplication. We need this factor for number theoretic reasons unlike case (1). It is known [7, p. 83–87] that this is a prehomogeneous vector space and the orbit of

(2.11) $$w = e_{234} + e_{567} + e_{125} + e_{136} + e_{147}$$

is Zariski open.

We use the operation $D_3$ defined in (1.14) again. For $x \in V$, we define

(2.12) $$S_x = x(\wedge, \otimes)D_3(x)(\wedge, \otimes)D_3(x) \in \wedge^7 W \otimes W \otimes W \cong W \otimes W,$$

where $(\wedge, \otimes)$ means the wedge product for the first factor and the tensor product for the second factor. Let $\phi : W \otimes W \to \mathrm{Sym}^2 W$ be the canonical map.

**Definition (2.13)** $Q_x = \phi(S_x)$.

We regard $Q_x$ and $x$ as elements of $\mathrm{Sym}^2(W^*)^* \cong \mathrm{Sym}^2 \mathrm{Im}(\widetilde{\mathbb{O}})^*$ and $\wedge^3(W^*)^*$ respectively. So $Q_x$ is a quadratic polynomial on $W^*$ and $x$ is an alternating trilinear form on $W^*$. We use the same notation $Q_x$ for the associated bilinear form on $W^*$. We define the action of $\mathrm{GL}(W)$ on $W^*$ by $gf(v) = f(g^{-1}v)$. In this way, we identify $\mathrm{GL}(W)$ and $\mathrm{GL}(W^*)$. Note that if we use bases for $W, W^*$ which are dual to each other and identify $\mathrm{GL}(W), \mathrm{GL}(W^*)$ with the set of $7 \times 7$ matrices, $g \in \mathrm{GL}(7)$ corresponds to ${}^tg^{-1}$ by this identification. It is easy to see that

$$(2.14) \qquad Q_{(t,g)x} = t^3 (\det g) g Q_x$$

(we are taking the determinant considering $g \in \mathrm{GL}(W)$).

We used the software "MAPLE" [1] to compute $Q_x$ for $x = w$ and

$$(2.15) \qquad w' = e_{234} + e_{346} + e_{127} - e_{145}.$$

For example, to compute $Q_{w'}$, we associate a differential form $dx_2 \wedge dx_3 \wedge dx_4 + \cdots$ to $w'$ and the input is as follows.

```
> with(difforms);
> defform(x1=0,x2=0,x3=0,x4=0,x5=0,x6=0,x7=0);
> w:= &^(d(x2),d(x3),d(x4))+&^(d(x3),d(x4),d(x6))
  +&^(d(x1),d(x2),d(x7))-&^(d(x1),d(x4),d(x5));
> v:= x2*d(x3)&^d(x4)-x3*d(x2)&^d(x4)+x4*d(x2)&^d(x3)
  +x3*d(x4)&^d(x5)-x4*d(x3)&^d(x5)+x5*d(x3)&^d(x4)
  +x1*d(x2)&^d(x7)-x2*d(x1)&^d(x7)+x7*d(x1)&^d(x2)
  -x1* d(x4)&^ d(x5)-x4*d(x1)&^d(x5)+x5*d(x1)&^d(x4);
> w&^v&^v;
```

The result is

$$(2.16) \qquad \begin{aligned} Q_w &= 6(-e_1^2 + e_2 e_5 + e_3 e_6 + e_4 e_7), \\ Q_{w'} &= 6(-e_1 e_4 + e_2 e_3). \end{aligned}$$

Since $Q_w$ is non-degenerate, the discriminant of $Q_x$ is a non-zero relative invariant polynomial of degree 21 and this reproves the existence of a relative invariant polynomial. Since $Q_w$ is irreducible, $Q_x$ is irreducible as a polynomial of $v \in W^*$. If $Q_x$ is divisible by a non-constant polynomial $p(x)$ of $x$, $p(x)$ is a relative invariant polynomial. Since $Q_{w'}$ is degenerate, $p(w') = 0$. But since $Q_{w'}$ is non-zero, this is a contradiction. So we get the following proposition.

**Proposition (2.17)** *As a polynomial of $(x, v) \in V \oplus W^*$, $Q_x(v)$ is irreducible.*

If $A$ is a normed $k$–algebra, $\mathrm{Re}(A)$ is contained in the center of $A$. So the structure of $A$ is determined by its restriction to $\mathrm{Im}(A)$. On $W^*$, we define a product structure $(\,\cdot\,)_x$ depending on $x$ by the equation

$$(2.18) \qquad 3x(v_1, v_2, v_3) = Q_x(v_1, (v_2 \cdot v_3)_x)$$

for all $v_1, v_2, v_3 \in W^*$. It is known [7, p. 86] that $(G, V)$ has a relative invariant polynomial $\Delta(x)$ of degree seven. Since $\mathrm{ch}\, k = 0$, we may assume that $\Delta(x) \in k[V]$ and $\Delta(w) = 6$. Clearly, $\Delta((t,g)x) = t^7(\det g)^3 \Delta(x)$. For $x \in V_k^{ss}$ and $v \in W_k^*$, we define

(2.19)
$$\|v\|_x = \Delta(x)^{-1} Q_x(v).$$

**Definition (2.20)** $\mathbb{O}_x$ *is the algebra* $k \oplus W^*$ *such that* $W^* = \text{Im}(\mathbb{O}_x)$. *The norm is defined by* $\|v\| = \|v\|_x$ *for* $v \in W^*$, *and the product* $v_1 v_2$ *for* $v_1, v_2 \in W^*$ *is defined by*
$$\text{Re}(v_1 v_2) = -\Delta(x)^{-1} Q_x(v_1, v_2), \ \text{Im}(v_1 v_2) = (v_1 \cdot v_2)_x.$$

If $x = w$, the product structure and the norm of $\mathbb{O}_w$ coincide with those of $\widetilde{\mathbb{O}}$ by Proposition (2.10). Therefore, $\mathbb{O}_w \cong \widetilde{\mathbb{O}}$ and is a normed $k$–algebra. Suppose $x, y \in V_k^{ss}$, $(t,g) \in G_{\bar{k}}$, and $y = (t,g)x$. Then we define $m_{x,y,(t,g)} : W^* \to W^*$ by

(2.21)
$$m_{x,y,(t,g)}(v) = t^2 (\det g) g v$$

for $v \in W^*$. Note that we are taking the determinant of $g$ considering $g \in \text{GL}(W)$.

**Proposition (2.22)** (1) *For all* $x \in V_k^{ss}$, $\mathbb{O}_x$ *is a normed $k$–algebra and is a $k$–form of* $\widetilde{\mathbb{O}}$.
(2) *If* $x, y \in V_k^{ss}$, $(t,g) \in G_{\bar{k}}$, *and* $y = (t,g)x$, $m_{x,y,(t,g)}$ *is an isomorphism of normed $\bar{k}$–algebras from* $\mathbb{O}_{x\bar{k}}$ *to* $\mathbb{O}_{y\bar{k}}$.

*Proof.* Let $m = m_{x,y,(t,g)}$. Then

$$\begin{aligned}
Q_y(gv_1, (m(v_2) \cdot m(v_3))_y) &= t^4 (\det g)^2 Q_y(gv_1, (gv_2 \cdot gv_3)_y) \\
&= 3t^4 (\det g)^2 y(gv_1, gv_2, gv_3) \\
&= 3t^4 (\det g)^2 g^{-1} y(v_1, v_2, v_3) \\
&= 3t^5 (\det g)^2 x(v_1, v_2, v_3) \\
&= t^5 (\det g)^2 Q_x(v_1, (v_2 \cdot v_3)_x) \\
&= t^2 \det g Q_y(gv_1, g(v_2 \cdot v_3)_x) \\
&= Q_y(gv_1, m(v_2 \cdot v_3)_x)
\end{aligned}$$

for all $v_1, v_2, v_3 \in W^*$.

Therefore, $(m(v_2) \cdot m(v_3))_y = m(v_2 \cdot v_3)_x$. Since

$$\begin{aligned}
\Delta(y)^{-1} Q_y(m(v_1), m(v_2)) &= (t^7 (\det g)^3)^{-1} t^4 (\det g)^2 \Delta(x)^{-1} Q_y(gv_1, gv_2) \\
&= (t^3 \det g)^{-1} \Delta(x)^{-1} t^3 \det g Q_x(v_1, v_2) \\
&= \Delta(x)^{-1} Q_x(v_1, v_2)
\end{aligned}$$

for all $v_1, v_2 \in W^*$, $m$ preserves the norm also.

Since $\mathbb{O}_w$ is a normed $k$–algebra, this proves both (1), (2). □

Let $\mathfrak{O}$ be the set of $k$–isomorphism classes of $k$–forms of $\widetilde{\mathbb{O}}$. We regard $\text{Aut}(\widetilde{\mathbb{O}}) \subset \text{GL}(\text{Im}(\widetilde{\mathbb{O}}))$.

**Lemma (2.23)** *If $g \in \mathrm{Aut}\,(\widetilde{\mathbb{O}})$, $\det g = 1$.*

*Proof.* Since $Q_w(gv_1, gv_2) = Q_w(v_1, v_2)$ for all $v_1, v_2 \in W^*$, $(\det g)^2 = 1$. Since $g(v_1 \cdot v_2)_w = (gv_1 \cdot gv_2)_w$, $gw = w$ by (2.18). It is proved in [7, p. 86] that if $g \in G_w \cap \mathrm{GL}(W)$, $(\det g)^3 = 1$. This implies $\det g = (\det g)^3/(\det g)^2 = 1$. □

**Proposition (2.24)** *The map*

$$\alpha_V : G_k \setminus V_k^{\mathrm{ss}} \ni x \to \mathbb{O}_x \in \mathfrak{D}$$

*is well defined and is bijective.*

*Proof.* If $x, y \in V_k^{\mathrm{ss}}$, $(t, g) \in G_k$, and $y = (t, g)x$, then $m_{x,y,(t,g)}$ is a $k$–isomorphism. Therefore, the above map is well defined. Suppose $x, y \in V_k^{\mathrm{ss}}$ and $g' : \mathbb{O}_x \to \mathbb{O}_y$ is a $k$–isomorphism. Regarding $g' \in \mathrm{GL}(W)_k$ (by considering the contragredient representation), let $t = (\det g')^3$, $g = (\det g')^{-1}g'$. Then $(t, g) \in G_k$ and $t^2 \det gg = g'$. Since $\Delta(x)^{-1}Q_x(v_1, v_2) = \Delta(y)^{-1}Q_y(g'v_1, g'v_2)$ for all $v_1, v_2 \in W^*$, by a similar argument as in Proposition (2.22), $Q_{(t,g)^{-1}y} = Q_x$. So

$$\begin{aligned}
3t^4(\det g)^2 g^{-1} y(v_1, v_2, v_3) &= 3y(gv_1, g'v_2, g'v_3) \\
&= Q_y(gv_1, (g'v_2 \cdot g'v_3)_y) \\
&= Q_y(gv_1, g'(v_2 \cdot v_3)_x) \\
&= t^2 \det g Q_y(gv_1, g(v_2 \cdot v_3)_x) \\
&= t^5(\det g)^2 Q_{(t,g)^{-1}}(v_1, (v_2 \cdot v_3)_x) \\
&= t^5(\det g)^2 Q_x(v_1, (v_2 \cdot v_3)_x) \quad \text{by the above remark} \\
&= 3t^5(\det g)^2 x(v_1, v_2, v_3)
\end{aligned}$$

for all $v_1, v_2, v_3 \in W^*$.

Therefore, $y = (t, g)x \in G_k x$.

Let $A$ be a $k$–form of $\widetilde{\mathbb{O}}$ corresponding to a cohomology class $c \in \mathrm{H}^1(k, \mathrm{Aut}\,(\widetilde{\mathbb{O}}))$. Suppose $c = \{h_\sigma\}_\sigma$. Then

$$\mathrm{Im}(A)_k = \{v \in \mathrm{Im}(\widetilde{\mathbb{O}})_{\bar{k}} \mid h_\sigma v^\sigma = v \text{ for all } \sigma \in \mathrm{Gal}(\bar{k}/k)\}.$$

Since $\mathrm{Aut}\,(\widetilde{\mathbb{O}}) \subset \mathrm{GL}(\mathrm{Im}(\widetilde{\mathbb{O}}))$ and $\mathrm{H}^1(k, \mathrm{GL}(\mathrm{Im}(\widetilde{\mathbb{O}})))$ is trivial, there exists $g' \in \mathrm{GL}(\widetilde{\mathbb{O}})_{\bar{k}}$ such that $c = \{g'^{-1}g'^\sigma\}_\sigma$. Then

$$\begin{aligned}
\mathrm{Im}(A)_k &= \{v \in \mathrm{Im}(\widetilde{\mathbb{O}})_{\bar{k}} \mid g'^{-1}g'^\sigma v^\sigma = v \text{ for all } \sigma \in \mathrm{Gal}(\bar{k}/k)\} \\
&= \{v \in \mathrm{Im}(\widetilde{\mathbb{O}})_{\bar{k}} \mid g'^\sigma v^\sigma = g'v \text{ for all } \sigma \in \mathrm{Gal}(\bar{k}/k)\} \\
&= \{v \in \mathrm{Im}(\widetilde{\mathbb{O}})_{\bar{k}} \mid g'v \in \mathrm{Im}(\widetilde{\mathbb{O}})_k\}.
\end{aligned}$$

So $A$ is characterized as the $k$–form such that $g'$ induces an isomorphism $\mathrm{Im}(\widetilde{\mathbb{O}})_{\bar{k}} \to \mathrm{Im}(A)_{\bar{k}}$.

Let $t = (\det g')^3$ and $g = (\det g')^{-1}g'$. Then $(t, g) \in G_{\bar{k}}$. By Lemma (2.23),
$$(t, g)^{-1}(t, g)^\sigma = (1, g'^{-1}g'^\sigma) \in G_{w\bar{k}}.$$

So this defines a cohomology class in $\mathrm{H}^1(k, G_w)$. This implies that $x = (t, g)w \in V_k^{\mathrm{ss}}$. Since $g' = m_{w,x,(t,g)}$, $g'$ induces an isomorphism from $\mathrm{Im}(\widetilde{\mathbb{O}})_{\bar{k}}$ to $\mathrm{Im}(\mathbb{O}_x)_{\bar{k}}$. Therefore, $A \cong \mathbb{O}_x$. □

**Remark (2.25)** It is proved in [7], [4] that $G_k \setminus V_k^{\mathrm{ss}} \cong \mathrm{H}^1(k, \mathrm{Aut}\,(\widetilde{\mathbb{O}}))$. So the credit for the existence of a bijective correspondence between $G_k \setminus V_k^{\mathrm{ss}}$ and $\mathfrak{O}$ should go to Sato–Kimura [7] and Igusa [4]. However, we constructed $\mathbb{O}_x \in \mathfrak{O}$ for $x \in V_k^{\mathrm{ss}}$, and the fact that this particular correspondence is bijective still required a proof. The operator $D_3$ was considered in [7]. The fact that the stabilizer of $w$ is a group of type $G_2$ at least goes back to [7].

For the rest of this section, we describe the orbit space $G_\mathbb{R} \setminus V_\mathbb{R}^{\mathrm{ss}} \cong \mathrm{GL}(W)_\mathbb{R} \setminus V_\mathbb{R}^{\mathrm{ss}}$. We will show that $G_\mathbb{R} \setminus V_\mathbb{R}^{\mathrm{ss}}$ consists of two orbits corresponding to $\mathbb{O}, \widetilde{\mathbb{O}}$.

It is known that $\mathbb{O}, \widetilde{\mathbb{O}}$ are the only $\mathbb{R}$–forms of $\widetilde{\mathbb{O}}$. We can identify $\mathbb{H}$ with $\mathbb{C}(+)$ by
$$a + bi + cj + dk \to (a + bi) + (c + di)\epsilon.$$

So $\mathbb{H} \otimes \mathbb{C}$ and $\mathrm{M}(2, 2)_\mathbb{C}$ are isomorphic by the map
$$a + bi + cj + dk \to \begin{pmatrix} a + \sqrt{-1}c & -b + \sqrt{-1}d \\ b + \sqrt{-1}d & a - \sqrt{-1}c \end{pmatrix}$$

for $a, b, c, d \in \mathbb{C}$ (see (2.6)).

Therefore, by considering $\mathbb{O} = \mathbb{H}(+)$, $\widetilde{\mathbb{O}} = \mathrm{M}(2, 2)(+)$, $\mathbb{O} \otimes \mathbb{C}$ and $\widetilde{\mathbb{O}} \otimes \mathbb{C}$ are isomorphic by the map

$$\phi((y_1 i + y_2 j + y_3 k) + (y_4 + y_5 i + y_6 j + y_7 k)\epsilon)$$
$$= \begin{pmatrix} \sqrt{-1}y_2 & -y_1 + \sqrt{-1}y_3 \\ y_1 + \sqrt{-1}y_3 & -\sqrt{-1}y_2 \end{pmatrix} + \begin{pmatrix} y_4 + \sqrt{-1}y_6 & -y_5 + \sqrt{-1}y_7 \\ y_5 + \sqrt{-1}y_7 & y_4 - \sqrt{-1}y_6 \end{pmatrix} \epsilon$$

Let $y = (y_1 i + y_2 j + y_3 k) + (y_4 + y_5 i + y_6 j + y_7 k)\epsilon$. We define

$$g_1 = \begin{pmatrix} 0 & -1 & 0 & 0 & -1 & 0 & 0 \\ \sqrt{-1} & 0 & 0 & 0 & 0 & 0 & 0 \\ 0 & \sqrt{-1} & 0 & 0 & -\sqrt{-1} & 0 & 0 \\ 0 & 0 & 1 & 0 & 0 & 1 & 0 \\ 0 & 0 & 0 & -1 & 0 & 0 & -1 \\ 0 & 0 & \sqrt{-1} & 0 & 0 & -\sqrt{-1} & 0 \\ 0 & 0 & 0 & -\sqrt{-1} & 0 & 0 & \sqrt{-1} \end{pmatrix}.$$

Consider $f_1, \cdots, f_7$ in (2.9). Suppose $x = \sum_{i=1}^7 x_i f_i = \phi(y)$. Then

$$\begin{pmatrix} y_1 \\ \vdots \\ y_7 \end{pmatrix} = {}^t g_1^{-1} \begin{pmatrix} x_1 \\ \vdots \\ x_7 \end{pmatrix}.$$

Note that $\det g_1 = 2^3$. So we put $t = 2^9$, $g = 2^{-3}g_1$ following the argument of Proposition (2.24). Then $t^2 \det gg = g_1$.

Consider $(t, g) \in \mathrm{GL}(W)$. Since $tgv = 2^6 g_1 v$ for $v \in W$, by $(t, g)$, $(e_1, \cdots, e_7)$ maps to $2^6(e_1, \cdots, e_7)g_1$. By an easy computation,

(2.26) $\qquad w_1 = (t, g)w = 2^7(e_{145} - e_{167} + e_{347} - e_{356} + e_{123} + e_{246} + e_{257}).$

By the proof of Proposition (2.24), $w_1$ corresponds to $\mathbb{O}$. It is easy to see that

(2.27) $\qquad\qquad\qquad \|\phi(y)\|_w = y_1^2 + \cdots + y_7^2.$

Since $\Delta(w_1) = 2^9 \cdot 6$, $Q_{w_1}(y) = 2^9 \cdot 6(y_1^2 + \cdots + y_7^2)$. Since $H^0_{w_1} \subset \mathrm{SO}(Q_{w_1})$, the real rank of $H^0_{w_1\mathbb{R}}$ is zero.

## §3 Intermediate groups

For $x \in V_k^{\mathrm{ss}}$, let $\mathfrak{h}_x$ be the Lie algebra of $H^0_x$. Consider the element $w$ for cases (1)–(3). Let $\mathfrak{h}_1 = \mathfrak{h}_w$ and $\mathfrak{h}_2 = \mathrm{sl}(6)$, $\mathrm{sl}(7)$, or $\mathrm{sl}(2n)$. In this section, we consider Lie subalgebras of $\mathfrak{h}_2$ containing $\mathfrak{h}_1$.

We first consider case (1). Clearly,

(3.1) $\qquad\qquad \mathfrak{h}_1 = \left\{ \begin{pmatrix} A & 0 \\ 0 & B \end{pmatrix} \,\Big|\, \mathrm{tr}(A) = \mathrm{tr}(B) = 0 \right\}.$

Let

(3.2) $\qquad\qquad \mathfrak{u}_1 = \left\{ \begin{pmatrix} 0 & U \\ 0 & 0 \end{pmatrix} \,\Big|\, U \in \mathrm{M}(3, 3) \right\},$

$\qquad\qquad \mathfrak{u}_2 = \left\{ \begin{pmatrix} 0 & 0 \\ U & 0 \end{pmatrix} \,\Big|\, U \in \mathrm{M}(3, 3) \right\},$

$\qquad\qquad \mathfrak{t} = \left\{ a \begin{pmatrix} I_3 & 0 \\ 0 & -I_3 \end{pmatrix} \,\Big|\, a \in k \right\}.$

Obviously $\mathfrak{h}_1$ is contained in the following Lie algebras

(3.3) $\qquad\qquad \mathfrak{h}_3 = \mathfrak{h}_1 \oplus \mathfrak{u}_1, \; \mathfrak{h}_4 = \mathfrak{h}_1 \oplus \mathfrak{u}_2,$
$\qquad\qquad \mathfrak{h}'_1 = \mathfrak{h}_1 \oplus \mathfrak{t}, \; \mathfrak{h}'_3 = \mathfrak{h}_1 \oplus \mathfrak{u}_1 \oplus \mathfrak{t}, \; \mathfrak{h}'_4 = \mathfrak{h}_1 \oplus \mathfrak{u}_2 \oplus \mathfrak{t}.$

**Proposition (3.4)** *If $\mathfrak{h}_1 \subset \mathfrak{f} \subset \mathfrak{h}_2$ is a Lie subalgebra, $\mathfrak{f} = \mathfrak{h}_1, \mathfrak{h}'_1, \mathfrak{h}_2, \mathfrak{h}_3, \mathfrak{h}'_3, \mathfrak{h}_4$, or $\mathfrak{h}'_4$.*

*Proof.* As a $\mathfrak{h}_1$–module, $\mathfrak{h}_2$ decomposes as a direct sum of representations as

(3.5) $\qquad\qquad\qquad \mathfrak{h}_2 = \mathfrak{h}_1 \oplus \mathfrak{u}_1 \oplus \mathfrak{u}_2 \oplus \mathfrak{t}.$

Let $\Lambda_1, \Lambda_2$ be the usual fundamental weights of $\mathrm{sl}(3)$, and $V_1, V_2$ the irreducible representations with highest weights $\Lambda_1, \Lambda_2$ respectively.

Then $\mathfrak{u}_1 \cong V_1 \otimes V_2$, $\mathfrak{u}_2 = V_2 \otimes V_1$, and $\mathfrak{t}$ is the trivial representation. Let $V_3$ be the representation of sl(3) with highest weight $\Lambda_1 + \Lambda_2$. Let $V_{3,1}$ be the representation of $\mathfrak{h}_1$ which is the tensor product of $V_3$ for the first factor of sl(3) and the trivial representation for the second factor of sl(3). We define $V_{3,2}$ similarly. Then $\mathfrak{h}_1$ is $V_{3,1} \oplus V_{3,2}$ as a representation of $\mathfrak{h}_1$. No two of these irreducible representations are equivalent. So

$$(3.6) \qquad \mathfrak{f} = \mathfrak{h}_1 \oplus (\mathfrak{f} \cap \mathfrak{u}_1) \oplus (\mathfrak{f} \cap \mathfrak{u}_2) \oplus (\mathfrak{f} \cap \mathfrak{t}).$$

Clearly, $\mathfrak{f} \cap \mathfrak{u}_1 = \mathfrak{u}_1$ or $0$, etc. Since $\mathfrak{h}_1, \mathfrak{u}_1, \mathfrak{u}_2$ generate $\mathfrak{h}_2$, if $\mathfrak{f}$ contains both $\mathfrak{u}_1, \mathfrak{u}_2$, $\mathfrak{f} = \mathfrak{h}_2$. This proves the proposition. □

We define $H_{w1} = H_w^0$, $H_{w2} = H$ and

$$(3.7) \qquad H_{w3} = \left\{ \begin{pmatrix} A & U \\ 0 & B \end{pmatrix} \middle| A, B \in \mathrm{SL}(3),\ U \in \mathrm{M}(3,3) \right\},$$

$$H_{w4} = \left\{ \begin{pmatrix} A & 0 \\ U & B \end{pmatrix} \middle| A, B \in \mathrm{SL}(3),\ U \in \mathrm{M}(3,3) \right\}.$$

If $k = \mathbb{R}$ and $x = gw$ for $g \in G_{\mathbb{R}}$, we define

$$(3.8) \qquad H_{x1\mathbb{R}} = H^0_{x\mathbb{R}+},\ H_{x2\mathbb{R}} = H_{\mathbb{R}},$$
$$H_{x3\mathbb{R}} = gH_{w3\mathbb{R}}g^{-1},\ H_{x4\mathbb{R}} = gH_{w4\mathbb{R}}g^{-1}.$$

This definition does not depend on the choice of $g \in G_{\mathbb{R}}$.

**Proposition (3.9)** *Suppose $x \in G_{\mathbb{R}} w$. Then if $H_{x1\mathbb{R}+} \subset F \subset H_{\mathbb{R}}$ is a closed connected subgroup whose radical is a unipotent subgroup, $F = H_{x1\mathbb{R}}, H_{x2\mathbb{R}}, H_{x3\mathbb{R}}$, or $H_{x4\mathbb{R}}$.*

*Proof.* Let $\mathfrak{f}$ be the Lie algebra of $F$. Then $\mathfrak{h}_{1\mathbb{R}} \subset g^{-1}\mathfrak{f}g \subset \mathfrak{h}_{2\mathbb{R}}$. So $g^{-1}\mathfrak{f}g = \mathfrak{h}_{1\mathbb{R}}, \mathfrak{h}'_{1\mathbb{R}}, \mathfrak{h}_{2\mathbb{R}}, \mathfrak{h}_{3\mathbb{R}}, \mathfrak{h}_{4\mathbb{R}}, \mathfrak{h}'_{3\mathbb{R}}$, or $\mathfrak{h}'_{4\mathbb{R}}$. But it cannot be $\mathfrak{h}'_{1\mathbb{R}}, \mathfrak{h}'_{3\mathbb{R}}$ or $\mathfrak{h}'_{4\mathbb{R}}$ because the radical of $F$ is a unipotent subgroup. This proves the proposition. □

We consider the orbit of $w_1$ (see (1.13)) next. Let $g_{w_1} = g_{\sqrt{-1}}$ (see (1.8)).

**Proposition (3.10)** *Let $x = gw_1$ for $g \in G_{\mathbb{R}}$. Then if $H^0_{x\mathbb{R}+} \subset F \subset H_{\mathbb{R}}$ is the connected component of the set of $\mathbb{R}$–rational points of an algebraic group defined over $\mathbb{Q}$ whose radical is a unipotent subgroup, $F = H^0_{x\mathbb{R}+}$ or $H_{\mathbb{R}}$.*

*Proof.* As in the previous proposition, we only have to consider the case $x = w_1$. Let $\mathfrak{f}$ be the Lie algebra of $F$, and $F_{\mathbb{C}} \subset H_{\mathbb{C}}$ the closed connected subgroup whose Lie algebra is $\mathfrak{f}_{\mathbb{C}}$. Then the radical of $F_{\mathbb{C}}$ is a unipotent subgroup and $F$ is the connected component of the set of $\mathbb{R}$–rational points of $F_{\mathbb{C}}$. So by a similar argument as in Proposition (3.9), $\mathfrak{f}_{\mathbb{C}} = \mathfrak{h}_{w_1 1 \mathbb{C}}, \mathfrak{h}_{w_1 2 \mathbb{C}}, \mathfrak{h}_{w_1 3 \mathbb{C}}$, or $\mathfrak{h}_{w_1 4 \mathbb{C}}$. We show that the last two cases cannot happen. Since the argument is similar, we only consider the case $\mathfrak{f} = \mathfrak{h}_{w_1 3 \mathbb{C}}$.

This implies that

$$F_{\mathbb{C}} = \left\{ g_{w_1} \begin{pmatrix} A & U \\ 0 & B \end{pmatrix} g_{w_1}^{-1} \middle| A, B \in \mathrm{SL}(3)_{\mathbb{C}},\ U \in \mathrm{M}(3,3)_{\mathbb{C}} \right\}.$$

The element $g_{w_1} \begin{pmatrix} A & U \\ 0 & B \end{pmatrix} g_{w_1}^{-1}$ is an $\mathbb{R}$–rational point if and only if $B = \overline{A}$ and $U = 0$. So
$$F = \left\{ g_{w_1} \begin{pmatrix} A & 0 \\ 0 & \overline{A} \end{pmatrix} g_{w_1}^{-1} \,\middle|\, A \in \mathrm{SL}(3)_{\mathbb{C}} \right\} = H^0_{w_1\mathbb{R}+}.$$
This implies $\mathfrak{f} = \mathfrak{h}_{w_1\mathbb{R}}$, which is a contradiction. $\square$

Next, we consider case (2). Let $\mathfrak{h}_3 = \mathrm{so}(Q_w)$.

**Proposition (3.11)** *Suppose $k$ is algebraically closed. Then if $\mathfrak{h}_1 \subset \mathfrak{f} \subset \mathfrak{h}_2$ is a Lie subalgebra, $\mathfrak{f} = \mathfrak{h}_1, \mathfrak{h}_2$, or $\mathfrak{h}_3$.*

*Proof.* Recall that $\mathfrak{h}_1$ is the simple Lie algebra of type $G_2$ and $\mathfrak{h}_1 \subset \mathfrak{h}_3 \subset \mathfrak{h}_2$. Let $W$ be the standard 7 dimensional representation of $\mathrm{sl}(7)$. As a representation of $\mathfrak{h}_1$, $W$ is irreducible and has highest weight equal to a fundamental weight $\Lambda_1$. There are no non-trivial representations of $\mathfrak{h}_1$ of lower dimension. The $\mathfrak{h}_1$–invariant complement to $\mathfrak{h}_1$ in $\mathfrak{h}_3$ is 7 dimensional and non-trivial (as $\mathfrak{h}_3$ is simple). Therefore $\mathfrak{h}_3 = \mathfrak{h}_1 \oplus W$. Now consider $\mathfrak{h}_3 \subset \mathfrak{h}_2$. Let $\theta$ be the order two automorphism of $\mathfrak{h}_2$ defined by $\theta(X) = -^*X$ (the adjoint with respect to $Q_w$). Then the fixed point set of $\theta$ is $\mathfrak{h}_3$. Write $\mathfrak{h}_2 = \mathfrak{h}_3 \oplus U$, the $\pm 1$ eigenspace decomposition for $\theta$. Then $[U, U]$ is contained in $\mathfrak{h}_3$ and is $\mathfrak{h}_3$–invariant, so $[U, U] = \mathfrak{h}_3$. Note that $\mathfrak{h}_2 \cong W \otimes W^*$ and $W \cong W^*$, as $\mathfrak{h}_2$–representations. Because the weight $2\Lambda_1$ does not occur in $\mathfrak{h}_3$, this implies that the representation $U$ has highest weight $2\Lambda_1$. The irreducible representation with highest weight $2\Lambda_1$ has dimension 27 [7, p. 21], the same as $U$, so $U$ is irreducible. Thus, $\mathfrak{h}_2 = \mathfrak{h}_1 \oplus W \oplus U$, as representations of $\mathfrak{h}_1$.

Now suppose that $\mathfrak{f}$ is a subalgebra with $\mathfrak{h}_1 \subset \mathfrak{f} \subset \mathfrak{h}_2$ with $\mathfrak{f} \neq \mathfrak{h}_3$. Then if $\mathfrak{f}$ contains $U$ it contains $\mathfrak{h}_3$, so must equal $\mathfrak{h}_2$. The proposition follows. $\square$

**Corollary (3.12)** *Suppose $x \in V_{\mathbb{R}}^{\mathrm{ss}}$. Then if $H^0_{x\mathbb{R}+} \subset F \subset H_{\mathbb{R}}$ is a closed connected subgroup, $F = H^0_{x\mathbb{R}+}, \mathrm{SO}(Q_x)_{\mathbb{R}}$, or $H_{\mathbb{R}}$.*

*Proof.* Let $\mathfrak{f}$ be the Lie algebra of $F$. Then $\mathfrak{h}_{x\mathbb{C}} \subset \mathfrak{f}_{\mathbb{C}} \subset \mathfrak{h}_{2\mathbb{C}}$ and $\mathfrak{f} = \mathfrak{f}_{\mathbb{C}} \cap \mathfrak{h}_2$. Suppose $x = gw$ where $g \in G_{\mathbb{C}}$. Then
$$g^{-1}\mathfrak{h}_{x\mathbb{C}}g = \mathfrak{h}_{1\mathbb{C}} \subset g^{-1}\mathfrak{f}_{\mathbb{C}}g \subset \mathfrak{h}_{2\mathbb{C}}.$$
So $g^{-1}\mathfrak{f}_{\mathbb{C}}g = \mathfrak{h}_{1\mathbb{C}}, \mathrm{so}(Q_w)_{\mathbb{C}}$, or $\mathfrak{h}_{2\mathbb{C}}$. Therefore,
$$\mathfrak{f}_{\mathbb{C}} = g\mathfrak{h}_{1\mathbb{C}}g^{-1} = \mathfrak{h}_{x\mathbb{C}}, g\mathrm{so}(Q_w)_{\mathbb{C}}g^{-1} = \mathrm{so}(Q_x)_{\mathbb{C}}, \text{ or } \mathfrak{h}_{2\mathbb{C}}.$$
This implies that
$$\mathfrak{f} = \mathfrak{h}_{x\mathbb{R}}, \mathrm{so}(Q_x)_{\mathbb{R}}, \text{ or } \mathfrak{h}_{2\mathbb{R}}.$$

$\square$

Finally, we consider case (3).

**Proposition (3.13)** *If $H^0_{w\mathbb{R}+} \subset F \subset H_{\mathbb{R}}$ is a closed connected subgroup, $F = H^0_{w\mathbb{R}+}$ or $H_{\mathbb{R}}$.*

*Proof.* Let $\mathfrak{f}$ be the Lie algebra of $F$. Then $\mathfrak{h}_1 \subset \mathfrak{f} \subset \mathfrak{h}_2$. Let $\Lambda_1, \Lambda_2$ be the first and the second fundamental weights of $\mathfrak{h}_1$. Then the standard representation $V_1 = W$ of $\mathfrak{h}_1$ has the highest weight $\Lambda_1$. Let $V_2$ be the irreducible representation of $\mathfrak{h}_1$ with highest weight $\Lambda_2$. Since $w$ is a $G_w$–invariant non-degenerate bilinear form on $W^*$, $W$ is equivalent to $W^*$. So $\mathfrak{h}_2$ plus the trivial representation is equivalent to $W \otimes W$. Clearly $W \otimes W \cong \wedge^2 W \oplus \operatorname{Sym}^2 W$. It is known (see [7, pp. 14,15]) that $\wedge^2 W$ is a sum of $V_2$ and the trivial representation. As a representation of $\mathfrak{h}_1$, $\mathfrak{h}_1$ is irreducible and has the same highest weight $2\Lambda_1$ as $\operatorname{Sym}^2 W$. By eliminating the trivial representation,
$$\mathfrak{h}_2 = \mathfrak{h}_1 \oplus V_2.$$
Since $\mathfrak{h}_1, V_2$ are not equivalent, $\mathfrak{f} = \mathfrak{h}_1 \oplus (\mathfrak{f} \cap V_2)$. So $\mathfrak{f} = \mathfrak{h}_1$ or $\mathfrak{h}_2$. $\square$

Since $V_{\mathbb{R}}^{\mathrm{ss}}$ is a single $G_{\mathbb{R}}$–orbit in this case, the following is an immediate consequence of the above proposition.

**Corollary (3.14)** *Suppose $x \in V_{\mathbb{R}}^{\mathrm{ss}}$. Then if $H_{x\mathbb{R}+}^0 \subset F \subset H_{\mathbb{R}}$ is a closed connected subgroup, $F = H_{x\mathbb{R}+}^0$, or $H_{\mathbb{R}}$.*

## §4 The fixed point set of $H_x^0$

We consider the fixed point set of $H_x^0$ for $x \in V_k^{\mathrm{ss}}$. Throughout this section, we assume that $k$ is algebraically closed.

We first consider case (1). Let $W$ be the standard representation of $\operatorname{GL}(6)$.

**Proposition (4.1)** *If $y \in V_k$ is fixed by $H_{wk}^0$, there exist $\alpha_1, \alpha_2 \in k$ such that $x = \alpha_1 e_{123} + \alpha_2 e_{456}$.*

*Proof.* We use the notation of (3.1)–(3.6). Let $V_{1,1}$ be the representation of $\mathfrak{h}_1$ which is the tensor product of $V_1$ for the first factor and the trivial representation for the second factor. We define $V_{1,2}$ similarly. Then $W = V_{1,1} \oplus V_{1,2}$. Therefore,

$$(4.2) \qquad \wedge^3 W \cong \wedge^3 V_{1,1} \oplus (\wedge^2 V_{1,1} \otimes V_{1,2}) \oplus (V_{1,1} \otimes \wedge^2 V_{1,2}) \oplus \wedge^3 V_{1,2}.$$

The second and the third factors are irreducible and non-trivial, and the first and the fourth factors are the trivial representation. Therefore, the dimension of the subspace
$$\{y \in V \mid gy = y \text{ for all } g \in H_{wk}^0\}$$
is two.

Obviously, $\alpha_1 e_{123} + \alpha_2 e_{456}$ belongs to the above space for all $\alpha_1, \alpha_2 \in k$. This proves the proposition. $\square$

The following is an immediate consequence of the above proposition.

**Corollary (4.3)** *Let $x \in V_k^{\mathrm{ss}}$. Then if $y \in V_k$ is fixed by $H_{xk}^0$, $\operatorname{Gr}(y) = \operatorname{Gr}(x)$.*

We consider case (2) next. Let $W$ be the standard representation of $\operatorname{GL}(7)$.

**Proposition (4.4)** *If $y \in V_k$ is fixed by $H_{wk}^0$, $y$ is a scalar multiple of $w$.*

*Proof.* Recall that $H_w^0$ is a simple group of type $G_2$ and $W$ is the 7 dimensional irreducible representation with highest weight $\Lambda_1$. We show that the 35 dimensional representation $\wedge^3 W$ contains the trivial representation $k$ exactly once. The weights of $W$ are the short roots of $\mathfrak{g}_2$, together with zero, so the non-zero weights form the vertices of a regular hexagon centered at 0. Therefore, by listing the non-zero weights $\alpha_1, \ldots, \alpha_6$ in the order they appear around the hexagon (with $\alpha_1 = \Lambda_1$), we have $\alpha_{i-1} + \alpha_{i+1} = \alpha_i$ for $1 \leq i \leq 6$ (with subscripts read modulo 6). So one sees that the highest weight of $\wedge^3 W$ is $2\Lambda_1$ (from the sum $\alpha_6 + \alpha_1 + \alpha_2$) and the multiplicity of the zero weight is 5 (from the three sums of the form $-\alpha_i + 0 + \alpha_i$ and the two sums of the form $\alpha_i + \alpha_{i+2} + \alpha_{i+4}$). Since the irreducible representation with highest weight $2\Lambda_1$ has dimension 27 and must occur, 8 dimensions remain. As the smallest non-trivial representation of $G_2$ has dimension 7, the only possibility is $W$ plus one copy of the trivial representation. Thus, $\wedge^3 W = U \oplus W \oplus k$. □

The following is an immediate consequence of the above proposition.

**Corollary (4.5)** *Let $x \in V_k^{\mathrm{ss}}$. Then if $y \in V_k$ is fixed by $H_{xk}^0$, $y$ is a scalar multiple of $x$.*

Finally we consider case (3). Let $W$ be the standard representation of $\mathrm{Sp}(2n)$. As we pointed out in §3, $\wedge^2 W$ is a sum of a non-trivial irreducible representation and the trivial representation. Therefore, the following proposition follows by an argument as above.

**Proposition (4.6)** *Let $x \in V_k^{\mathrm{ss}}$. Then if $y \in V_k$ is fixed by $H_{xk}^0$, $y$ is a scalar multiple of $x$.*

### §5 The orbit closure

In this section, we formulate irrationality conditions for $x \in V_{\mathbb{R}}^{\mathrm{ss}}$ for cases (1)–(3) and prove that the orbit of $H_{x\mathbb{R}+}^0$ in $H_{\mathbb{R}}/\Gamma$ is dense if $x \in V_{\mathbb{R}}^{\mathrm{ss}}$ is sufficiently irrational.

Consider case (1). Let $G, V, W, w_0 = w, w_1$, etc. be as in §1.

**Definition (5.1)(1)** *A point $x \in G_{\mathbb{R}} w$ is sufficiently irrational if both $[E_{x1}], [E_{x2}] \in \mathrm{Gr}(3,6)_{\mathbb{C}}$ are irrational and $\mathrm{Gr}(x) \in ((\mathbb{Z}/2\mathbb{Z}) \backslash (\mathrm{Gr}(3,6) \times \mathrm{Gr}(3,6)))_{\mathbb{R}}$ is irrational. A point $x \in G_{\mathbb{R}} w_1$ is sufficiently irrational if $\mathrm{Gr}(x) \in ((\mathbb{Z}/2\mathbb{Z}) \backslash (\mathrm{Gr}(3,6) \times \mathrm{Gr}(3,6)))_{\mathbb{R}}$ is irrational.*

Note that if $x \in G_{\mathbb{R}} w_1$, $E_{x1}, E_{x2}$ are complex conjugates of each other. So the irrationality of one of $[E_{x1}], [E_{x2}]$ implies the irrationality of the other. Also even though we cannot distinguish $[E_{x1}], [E_{x2}]$, if we say both are irrational, the statement makes sense.

Next consider case (2). Let $G, V, W$, etc. be as in §2.

**Definition (5.1)(2)** *A point $x \in V_{\mathbb{R}}^{\mathrm{ss}}$ is sufficiently irrational if $[Q_x] \in \mathbb{P}(\mathrm{Sym}^2 W)_{\mathbb{R}}$ is irrational.*

Next consider case (3). Let $W$ be the standard representation of $\mathrm{GL}(2n)$ and $V = \wedge^2 W$.

**Definition (5.1)(3)** *A point $x \in V_{\mathbb{R}}^{\mathrm{ss}}$ is sufficiently irrational if $[x] \in \mathbb{P}(V)_{\mathbb{R}}$ is irrational.*

Now we are ready to prove that the orbit of $H_{x\mathbb{R}+}^0$ is dense.

**Theorem (5.2)** *Let $\Gamma \subset \mathrm{SL}(W)_\mathbb{R}$ be an arithmetic lattice. Then if $x \in V_\mathbb{R}^{ss}$ is sufficiently irrational in the sense of* Definition (5.2)(1)–(3) *and the real rank of $H^0_{x\mathbb{R}}$ is positive, $H^0_{x\mathbb{R}+}\Gamma$ is dense in $H_\mathbb{R}$.*

*Proof.* Note that the real rank of $H^0_{x\mathbb{R}}$ is positive for all $x \in V_\mathbb{R}^{ss}$ for cases (1), (3), and for $x \in G_\mathbb{R} w$ for case (2).

The proof for case (2) is the same as the proof of Theorem (5.1) [10] (using Corollaries (3.12),(4.5)). The proof for case (3) is also similar using Corollary (3.14) and Proposition (4.6). So we only consider case (1).

By Ratner's theorem (see the introduction of [10]), there exists a connected closed subgroup $H^0_{x\mathbb{R}+} \subset F \subset \mathrm{SL}(W)_\mathbb{R}$ such that $\overline{H^0_{x\mathbb{R}+}\Gamma} = F\Gamma$. By Proposition (3.2) [8, pp. 321–322], $F$ is defined over $\mathbb{Q}$ and the radical of $F$ is a unipotent subgroup. So by Propositions (3.9), (3.10), $F = H_{x1\mathbb{R}}, H_{x2\mathbb{R}}, H_{x3\mathbb{R}}$, or $H_{x4\mathbb{R}}$.

Suppose $F = H^0_{x\mathbb{R}+}$. Since $F$ is defined over $\mathbb{Q}$, for any $\sigma \in \mathrm{Aut}\,(\mathbb{C}/\mathbb{Q})$, $F^\sigma_\mathbb{C} = F_\mathbb{C}$. So $H^0_{x^\sigma \mathbb{C}} = H^0_{x\mathbb{C}} = H^0_{x\mathbb{C}}$. Since this group fixes $x^\sigma$,
$$\mathrm{Gr}(x)^\sigma = \mathrm{Gr}(x^\sigma) = \mathrm{Gr}(x)$$
by Proposition (1.20) and Corollary (4.3). Since this is the case for all $\sigma$, $\mathrm{Gr}(x)$ is a $\mathbb{Q}$–rational point of $(\mathbb{Z}/2\mathbb{Z}) \setminus (\mathrm{Gr}(3,6) \times \mathrm{Gr}(3,6))$, which is a contradiction.

Suppose $F = H_{x3\mathbb{R}}$ or $H_{x4\mathbb{R}}$ (which means $x \in G_\mathbb{R} w$). Since the argument is similar, we only consider the case $F = H_{x3\mathbb{R}}$. Suppose $x = gw$ for $g \in \mathrm{GL}(W)_\mathbb{C}$. Then if $\sigma \in \mathrm{Aut}\,(\mathbb{C}/\mathbb{Q})$, $x^\sigma = g^\sigma w$ and
$$\begin{aligned}
H^\sigma_{x3\mathbb{C}} &= \left\{ g^\sigma \begin{pmatrix} A^\sigma & U^\sigma \\ 0 & B^\sigma \end{pmatrix} (g^\sigma)^{-1} \,\bigg|\, A, B \in \mathrm{SL}(3)_\mathbb{C},\ U \in \mathrm{M}(3,3)_\mathbb{C} \right\} \\
&= \left\{ g^\sigma \begin{pmatrix} A & U \\ 0 & B \end{pmatrix} (g^\sigma)^{-1} \,\bigg|\, A, B \in \mathrm{SL}(3)_\mathbb{C},\ U \in \mathrm{M}(3,3)_\mathbb{C} \right\} \\
&= H_{x^\sigma 3\mathbb{C}}.
\end{aligned}$$

Since this group is defined over $\mathbb{Q}$, $H_{x3\mathbb{C}} = H^\sigma_{x3\mathbb{C}} = H_{x^\sigma 3\mathbb{C}}$. Note that $[E_{w1}] \in \mathrm{Gr}(3,6)_\mathbb{C}$ is the unique point fixed by $H_{w3\mathbb{C}}$. Since $g[E_{w1}]$ (resp. $g^\sigma[E_{w1}]$) is fixed by $H_{x3\mathbb{C}}$ (resp. $H_{x^\sigma 3\mathbb{C}}$), $g[E_{w1}] = g^\sigma[E_{w1}] = (g[E_{w1}])^\sigma$. Since this is the case for all $\sigma$, $g[E_{w1}]$ must be a $\mathbb{Q}$–rational point of $\mathrm{Gr}(3,6)$. So either $[E_{x1}]$ or $[E_{x2}]$ is a $\mathbb{Q}$–rational point of $\mathrm{Gr}(3,6)$, which is a contradiction. This proves that $F = H_\mathbb{R}$. □

## §6 Analogues of the Oppenheim conjecture

We prove our main theorem of this paper in this section.

We consider case (1) in (6.1)–(6.6). Consider $V, W, \{e_1, \cdots, e_6\}, v, w_0 = w, w_1$ in §1.

**Proposition (6.1)** *Let $y = (y_{ijk}) \in \wedge^3 \mathbb{R}^5$, and $\epsilon > 0$. Then there exist $z_0 = (z_{0,ijk}) \in G_\mathbb{R} w_0$ and $z_1 = (z_{1,ijk}) \in G_\mathbb{R} w_1$ such that $|y_{ijk} - z_{l,ijk}| < \epsilon$ for $l = 0, 1,\ 1 \leq i < j < k \leq 5$.*

*Proof.* Let $z = (z_{ijk})$ and

(6.2) $$X = \begin{pmatrix} z_{234} & -z_{134} & z_{124} \\ z_{235} & -z_{135} & z_{125} \\ z_{236} & -z_{136} & z_{126} \end{pmatrix},\ Y = \begin{pmatrix} z_{156} & -z_{146} & z_{145} \\ z_{256} & -z_{246} & z_{245} \\ z_{356} & -z_{346} & z_{345} \end{pmatrix}.$$

Note that $z \in G_\mathbb{R} w_0$ (resp. $z \in G_\mathbb{R} w_1$) if and only if $\Delta(z) > 0$ (resp. $\Delta(z) < 0$). It is proved in [7, p. 83] that

$$(6.3) \quad \Delta(z) = (z_{123} z_{456} - \mathrm{tr}(XY))^2 + 4z_{123} \det Y \\ + 4z_{123} \det X - 4 \sum_{i,j} \det X_{ij} \det Y_{ji},$$

where $X_{ij}, Y_{ji}$ are the $(i,j)$–minor and the $(j,i)$–minor of $X, Y$ respectively.

As a polynomial of $z_{456}$, this is a quadratic polynomial and the coefficient of $z_{456}^2$ is $z_{123}^2$. Given $y$, we can choose $z_{123} \neq 0$ close to $y_{123}$. By taking $z_{456} \gg 0$, $\Delta(z) > 0$. So the existence of $z_0$ follows.

We try to choose $z_1$ so that $z_{1,ij6} = 0$ unless $(i,j) = (1,5), (2,4), (4,5)$. So we assume $z_{ij6} = 0$ unless $(i,j) = (1,5), (2,4), (4,5)$ in the following. We first choose $z_{123} \neq 0$ close to $y_{123}$. Then

$$(6.4) \quad \Delta(z) = z_{123}^2 (z_{456} - z_{123}^{-1} \mathrm{tr}(XY))^2 + 4z_{123} \det Y \\ + 4(z_{456} - z_{123}^{-1} \mathrm{tr}(XY)) \det X \\ + 4z_{123}^{-1} \mathrm{tr}(XY) \det X - 4 \sum_{i,j} \det X_{ij} \det Y_{ji}.$$

As a quadratic polynomial of $z_{456} - z_{123}^{-1} \mathrm{tr}(XY)$, (6.4) can take a negative value if the discriminant is positive. Let

$$(6.5) \quad f_1(z) = 16 z_{123}^2 (-z_{123} z_{345} + z_{234} z_{135} + z_{145} z_{235})$$

Then a simple consideration shows that there exist polynomials $f_2(z), f_3(z), f_4(z)$ which depend only on $(z_{ijk})_{i<j<k\leq 5}$ such that the the discriminant of (6.4) is

$$(6.6) \quad f_1(z) z_{156} z_{246} + f_2(z) z_{156} + f_3(z) z_{246} + f_4(z).$$

Note that if $f$ is a non-zero polynomial on a real vector space $U_\mathbb{R}$, the set $\{p \in U_\mathbb{R} \mid f(p) \neq 0\}$ is open dense in $U_\mathbb{R}$ in classical topology. Since $f_1(z)$ is not identically zero, we can choose $(z_{ijk})_{i<j<k\leq 5}$ close to $y$ so that $f_1(z) \neq 0$. If $f_1(z) > 0$, we choose $z_{156} = z_{246} \gg 0$. If $f_1(z) < 0$, we choose $z_{156} = -z_{246} \gg 0$. In both cases, we can make (6.6) positive. $\square$

We consider case (2) in (6.7)–(6.13). Consider $V, W, \{e_1, \cdots, e_7\}, v, w$ in §2. We use the notation $e_{i_1 \ldots i_k} = e_{i_1} \wedge \cdots \wedge e_{i_k}$ as before.

**Proposition (6.7)** *Let $y = (y_{ijk}) \in \wedge^3 \mathbb{R}^6$, and $\epsilon > 0$. Then there exists $z = (z_{ijk}) \in G_\mathbb{R} w$ such that $|y_{ijk} - z_{ijk}| < \epsilon$ for $1 \leq i < j < k \leq 6$.*

*Proof.* We try to choose $z$ so that $z \in V_\mathbb{R}^{\mathrm{ss}}$, $Q_z(v_1) > 0$, and $Q_z(v_7) < 0$. This condition ensures that $Q_z$ is indefinite and therefore, $z \in G_\mathbb{R} w$. Let

$$(6.8) \quad I_1 = \left\{ (j, k, j', k', j'', k'') \,\middle|\, \begin{array}{c} j < k, j' < k', j'' < k'', \\ \{j, k, j', k', j'', k''\} = \{2, \cdots, 7\} \end{array} \right\},$$

$$I_2 = \left\{ (i, j, i', j', i'', j'') \,\middle|\, \begin{array}{c} i < j, i' < j', i'' < j'', \\ \{i, j, i', j', i'', j''\} = \{1, \cdots, 6\} \end{array} \right\},$$

$$I_3 = \{(j, k, j'k') \mid j < k, j' < k', \{j, k, j', k'\} = \{3, 4, 5, 6\}\}.$$

**Lemma (6.9)** (1) *The coefficient of $v_1^2$ in $Q_z$ is*

$$f_1(z) = \sum_{(j,k,j',k',j'',k'') \in I_1} \operatorname{sgn}\begin{pmatrix} 2 & \cdots & 7 \\ j & \cdots & k'' \end{pmatrix} z_{1jk} z_{1j'k'} z_{1j''k''}.$$

(2) *The coefficient of $v_7^2$ in $Q_z$ is*

$$f_2(z) = \sum_{(i,j,i',j',i'',j'') \in I_2} \operatorname{sgn}\begin{pmatrix} 1 & \cdots & 6 \\ i & \cdots & j'' \end{pmatrix} z_{ij7} z_{i'j'7} z_{i''j''7}.$$

*Proof.* We consider (1) first. Consider $z \wedge D_3(z) \wedge D_3(z)$.

(6.10) $$D_3(z) = \sum_{i'<j'<k'} z_{i'j'k'}(e_{j'k'} \otimes e_{i'} - e_{i'k'} \otimes e_{j'} + e_{i'j'} \otimes e_{k'}).$$

Since $j', k' > 1$, the terms in the right hand side of (6.10) which has $e_1$ as the second factor are those with $i' = 1$. Therefore,

(6.11) $\quad z \wedge D_3(z) \wedge D_3(z)$

$$= \left( \sum_{i,j,k} \sum_{j',k'} \sum_{j'',k''} x_{ijk} x_{1j'k'} x_{1j''k''} e_{ijk} \wedge e_{j'k'} \wedge e_{j''k''} \right) \otimes e_1 \otimes e_1$$
$+ \cdots.$

Clearly, $e_{ijk} \wedge e_{j'k'} \wedge e_{j''k''} = 0$ unless $i = 1$ and $(j,k,j',k',j'',k'') \in I_1$. So we only have to consider indices in $I_1$. If $(j,k,j',k',j'',k'') \in I_1$,

$$e_{ijk} \wedge e_{j'k'} \wedge e_{j''k''} = \operatorname{sgn}\begin{pmatrix} 2 & \cdots & 7 \\ j & \cdots & k'' \end{pmatrix} e_{234567}.$$

This proves (1).

Note that

$$e_{ij7} \wedge e_{i'j'} \wedge e_{i''j''} = e_{ij} \wedge e_{i'j'} \wedge e_{i''j''} \wedge e_7.$$

So (2) is similar. □

We consider $z$ such that $z_{ij7} = 0$ unless $(i,j) = (1,2), (3,4), (5,6)$. Note that this has no effect on $z_{ijk}$ with $i < j < k \leq 6$. We define

(6.12) $$f_3(z) = \sum_{(j,k,j',k') \in I_3} \operatorname{sgn}\begin{pmatrix} 3 & \cdots & 6 \\ j & \cdots & k' \end{pmatrix} z_{1jk} z_{1j'k'}.$$

The following is an immediate consequence of the above Lemma.

**Corollary (6.13)** *Suppose $z_{ij7} = 0$ unless $(i,j) = (1,2), (3,4), (5,6)$. Then there exists a polynomial $f_4(z)$ which does not depend on $z_{127}, z_{347}, z_{567}$ such that*

$$f_1(z) = 3z_{127} f_3(z) + f_4(z), \quad f_2(z) = 6 z_{127} z_{347} z_{567}.$$

Since $f_3(z)$ is non-zero, we can choose $(z_{ijk})_{i<j<k\leq 6}$ arbitrarily close to $y$ and $f_3(z) \neq 0$. If $f_3(z) > 0$, we choose $z_{127} \gg 0$, $z_{347}z_{567} < 0$. If $f_3(z) > 0$, we choose $z_{127} \ll 0$, $z_{347}z_{567} > 0$. In both cases, $f_1(z) > 0, f_2(z) < 0$. Since $V_\mathbb{R}^{ss}$ is open dense in $V_\mathbb{R}$, we can replace $z$ if necessary and assume that $z \in V_\mathbb{R}^{ss}$. In this process, the condition $f_1(z) > 0, f_2(z) < 0$ can be preserved. This completes the proof of Proposition (6.7). □

Now we consider case (3). So $G = \mathrm{GL}(2n)$, $W = \mathbb{Q}^{2n}$, $V = \wedge^2 W$. In this case, $V_\mathbb{R}^{ss}$ is a single $G_\mathbb{R}$–orbit. So the following proposition is obvious.

**Proposition (6.14)** *Let $y = (y_{ij}) \in \wedge^2 \mathbb{R}^{2n-1}$, and $\epsilon > 0$. Then there exist $z = (z_{ij}) \in V_\mathbb{R}^{ss}$ such that $|y_{ij} - z_{ij}| < \epsilon$ for $1 \leq i < j \leq 2n-1$.*

Now we are ready to prove our main theorem. In statements (1)–(3) of the following theorem, we consider the prehomogeneous vector spaces (1)–(3) of this paper respectively. Let $W$ be the standard representation of $\mathrm{GL}(6), \mathrm{GL}(7)$, or $\mathrm{GL}(2n)$. We identify $\wedge^3 W$ and $\wedge^2 W$ with $\wedge^3(W^*)^*$ and $\wedge^2(W^*)^*$ respectively.

**Theorem (6.15)** *Suppose $x \in V_\mathbb{R}^{ss}$ is sufficiently irrational in the sense of Definition (5.1)(1)–(3) and the real rank of $H_{x\mathbb{R}}^0$ is positive.*
*(1) For any $y = (y_{ijk}) \in \wedge^3 \mathbb{R}^5$ and $\epsilon > 0$, there exists a $\mathbb{Z}$–basis $\{u_1, \cdots, u_6\}$ of $W_\mathbb{Z}^* \cong \mathbb{Z}^6$ such that*
$$|y_{ijk} - x(u_i, u_j, u_k)| < \epsilon$$
*for all $i < j < k \leq 5$.*
*(2) For any $y = (y_{ijk}) \in \wedge^3 \mathbb{R}^6$ and $\epsilon > 0$, there exists a $\mathbb{Z}$–basis $\{u_1, \cdots, u_7\}$ of $W_\mathbb{Z}^* \cong \mathbb{Z}^7$ such that*
$$|y_{ijk} - x(u_i, u_j, u_k)| < \epsilon$$
*for all $i < j < k \leq 6$.*
*(3) For any $y = (y_{ij}) \in \wedge^2 \mathbb{R}^{2n-1}$ and $\epsilon > 0$, there exists a $\mathbb{Z}$–basis $\{u_1, \cdots, u_{2n}\}$ of $W_\mathbb{Z}^* \cong \mathbb{Z}^{2n}$ such that*
$$|y_{ij} - x(u_i, u_j)| < \epsilon$$
*for all $i < j \leq 2n-1$.*

*Proof.* Since the proof is similar we only consider case (1). For cases (2), (3) the argument is similar using Theorem (5.2) and Propositions (6.7), (6.14).

By Proposition (6.1), we can choose $z = (z_{ijk}) \in G_\mathbb{R} x$ such that $|y_{ijk} - z_{ijk}| < \frac{\epsilon}{2}$ for $1 \leq i < j < k \leq 5$. Since any element of $G_\mathbb{R}$ can be written as a product of an element of $H_\mathbb{R}$ and a scalar matrix, we can choose $h' \in H_\mathbb{R}$ and $\lambda \in \mathbb{R} \setminus \{0\}$ so that $h'^{-1}x = \lambda z$. Let

$$h = h' \begin{pmatrix} \lambda^{\frac{1}{3}} & & & \\ & \ddots & & \\ & & \lambda^{\frac{1}{3}} & \\ & & & \lambda^{-\frac{5}{3}} \end{pmatrix} \in H_\mathbb{R}.$$

We put $z' = h^{-1}x$. Then if $z' = (z'_{ijk})$, $z'_{ijk} = z_{ijk}$ for all $i < j < k \leq 5$.

By Theorem (5.2), $H_{x\mathbb{R}+}^0 H_\mathbb{Z}$ is dense in $H_\mathbb{R}$. So we can choose $h_1 \in H_{x\mathbb{R}+}^0$ and $h_2 \in H_\mathbb{Z}$ so that $h_1 h_2$ is close to $h$. Then $(h_1 h_2)^{-1} x = h_2^{-1} h_1^{-1} x = h_2^{-1} x$ is close to $z' = h^{-1}x$. Let $f_1, \cdots, f_6$ be the standard coordinate vectors of $W^*$. Since $h_2^{-1}x$ is close to $z'$, we can assume that

$$|z'(f_i, f_j, f_k) - h_2^{-1}x(f_i, f_j, f_k)| = |z_{ijk} - x(h_2 f_i, h_2 f_j, h_2 f_k)| < \frac{\epsilon}{2}$$

for all $i < j < k \leq 5$.

Let $u_1 = h_2 f_1, \cdots, u_6 = h_2 f_6$. Since $h_2 \in H_{\mathbb{Z}}$, $\{u_1, \cdots, u_6\}$ is a $\mathbb{Z}$–basis of $W_{\mathbb{Z}}^*$. Moreover,

$$|y_{ijk} - x(u_i, u_j, u_k)| \leq |y_{ijk} - z_{ijk}| + |z_{ijk} - x(u_i, u_j, u_k)| < \epsilon$$

for all $i < j < k \leq 5$. This proves the theorem. $\square$

Dave Witte
Oklahoma State University
Mathematics Department
401 Math Science
Stillwater OK 74078-1058 USA
dwitte@math.okstate.edu

Akihiko Yukie
Oklahoma State University
Mathematics Department
401 Math Science

Stillwater OK 74078-1058 USA
yukie@math.okstate.edu
http://www.math.okstate.edu/˜yukie

Roger Zierau
Oklahoma State University
Mathematics Department
401 Math Science
Stillwater OK 74078-1058 USA
zierau@math.okstate.edu